\documentclass[a4paper,11pt,reqno,twoside]{article}
\usepackage{amssymb}
\usepackage[frenchb,english]{babel}
\usepackage{amscd}
\usepackage{amsmath,amsthm,amsfonts,amssymb,graphicx}
\usepackage[colorlinks,linkcolor=blue,anchorcolor=blue,citecolor=green]{hyperref}
\usepackage{mathrsfs}
\usepackage{fancyhdr}
\usepackage{multicol}
\usepackage{abstract}
\usepackage{mathrsfs}

\begin{document}
\allowdisplaybreaks[4]
\theoremstyle{plain}
\newtheorem{Definition}{Definition}[section]
\newtheorem{Proposition}{Proposition}[section]
\newtheorem{Property}{Property}[section]
\newtheorem{Theorem}{Theorem}[section]
\newtheorem{Lemma}{Lemma}[section]
\newtheorem{Corollary}[Theorem]{Corollary}
\newtheorem{Remark}{Remark}[section]

\noindent  {\LARGE
Special Toeplitz operators on a class of bounded Hartogs domains}\\\\
\noindent\text{Yanyan Tang \; \& \; Zhenhan Tu$^{*}$ }\\

\noindent\small {School of Mathematics and Statistics, Wuhan
University, Wuhan, Hubei 430072, P.R. China}

\noindent\text{Email: yanyantang@whu.edu.cn (Y. Tang), zhhtu.math@whu.edu.cn (Z. Tu)}

\renewcommand{\thefootnote}{{}}
\footnote{\hskip -16pt {$^{*}$Corresponding author. \\}}\\

\normalsize \noindent\textbf{Abstract}\quad  We introduce a wider class of bounded Hartogs domains,
which contains some generalizations of the classical Hartogs triangle.  A sharp criteria for the
$L^p-L^q$ boundedness  of the Toeplitz operator with symbol $K^{-t}$ is obtained on these domains,
where $K$ is the Bergman kernel on diagonal and $t\geq 0$. It generalizes the results by Chen and Beberok in the case $1<p<\infty$.

\vskip 10pt

\noindent \textbf{Key words:}  Hartogs domains  \textperiodcentered \; $L^p-L^q$ boundedness \textperiodcentered \; Toeplitz operator
\vskip 10pt

\noindent \textbf{Mathematics Subject Classification (2010):} 32A36  \textperiodcentered \, 32A25  \textperiodcentered \, 32A07

\pagenumbering{arabic}
\renewcommand{\theequation}
{\arabic{section}.\arabic{equation}}

\setcounter{section}{0}
\setcounter{equation}{0}
\section{Introduction}

\noindent 1.1. Toeplitz operator
\vspace{3mm}

Let $\Omega$ be a domain in $\mathbb{C}^n$ and
$A^2(\Omega)$ be the Bergman space of holomorphic functions in the square integrable space $L^2(\Omega)$. It is well known that the Bergman projection
$$
P_{\Omega}:  L^2(\Omega)\rightarrow A^2(\Omega)
$$
is an integral operator represented by
$$
P_{\Omega}(f)(z)=\int_{\Omega}K_{\Omega}(z,w)f(w)dA(w),
$$
where $K_{\Omega}(\cdot,\cdot)$ is the Bergman kernel on $\Omega\times\Omega$, and $dA$ is the Lebesgue measure on $\Omega$. The study of the $L^p$ boundedness of the  Bergman projection is a fact of interesting and fundamental importance. One of the most common object is the bounded domain with various boundary conditions.
For example, if $\Omega$ is a bounded strongly pseudoconvex domain or a smoothly bounded pseudoconvex domain of finite type in $\mathbb{C}^2$ or some bounded Reinhardt domains, then the Bergman projection $P_{\Omega}$ is bounded from $L^p(\Omega)$ to itself for all $p\in (1, \infty)$ (see Huo \cite{Huo}, Lanzani-Stein \cite{Lanzani-Stein}, McNeal \cite{McNeal}, and Phong-Stein \cite{Phong-Stein}). When the domain $\Omega$ has serious singularities at the boundary, in general, the $L^p$ boundedness of $P_{\Omega}$ will no longer hold for all $p\in (1, \infty)$ (see, e.g., Chakrabarti-Zeytuncu \cite{Chakrabarti-Zeytuncu}, Edholm-McNeal \cite{Edholm-McNeal}, Krantz-Peloso \cite{Krantz-Peloso}).
There are also smoothly bounded pseudoconvex domains where the Bergman projection has a restriction on $p$ for $L^p$ boundedness (see Barrett-\c{S}ahuto\u{g}lu \cite{Barrett}).

Given a function $\psi\in L^{\infty}(\Omega)$, the Toeplitz operator with symbol $\psi$ is defined by
$$
T_{\psi}(f)(z):=P_{\Omega}(\psi f)(z)=\int_{\Omega}K_{\Omega}(z,w)f(w)\psi(w)dA(w).
$$
When $\Omega$ is a bounded domain and $\psi$ is constant 1, then the Toeplitz operator will degenerate to the Bergman projection. Moreover, it is easy to see that $T_\psi: L^p(\Omega)\rightarrow L^q(\Omega)$ is bounded if $P_{\Omega}:L^p(\Omega)\rightarrow L^q(\Omega)$ is bounded for $\psi\in L^{\infty}(\Omega)$, where $p, q\in (1, \infty)$. It is natural to consider the following problem.

\vspace{2mm}
\noindent\textbf{Question} How does the symbol $\psi$ affect the boundedness of the Toeplitz operator $T_{\psi}$ from $L^p(\Omega)$ to $L^q(\Omega)$? Can we obtain a larger range of $p$ for the boundedness of $T_{\psi}$ in comparison with the corresponding Bergman projection?
\vspace{2mm}

The above question is the so-called ``gain" $L^p-L^q$ estimate properties of the Toeplitz operator.
In \u{C}u\u{c}kovi\'{c} and McNeal \cite{Cuckovic-McNeal}, they gave an affirmative answer to the above question on smoothly bounded strongly pseudoconvex domain $\Omega$ in $\mathbb{C}^n$ by choosing a special symbol $\psi:=\delta^\eta$, where $\eta\geq0$ and
$\delta(\cdot) £º=d(\cdot, \partial\Omega)$ is the Euclidean distance from the boundary. Later, Khanh-Liu-Thuc \cite{Khanh-Liu-Thuc1} also considered the same problem on certain classes of smoothly bounded pseudoconvex domains of finite type with symbol $\psi:=K^{-\alpha}$, where $\alpha\geq0$ and $K$
is the Bergman kernel on diagonal. Recently, they continued to work on the same question
on fat Hartogs triangle $\Omega_{k}=\{(z_1, z_2)\in \mathbb{C}^2: |z_1|^k<|z_2|<1\}~(k\in \mathbb{Z}^{+})$, which has singular boundary points
(see Khanh-Liu-Thuc \cite{Khanh-Liu-Thuc2}).

\vspace{3mm}
\noindent 1.2. Generalizations of the Hartogs triangle
\vspace{3mm}

Besides the power-generalized Hartogs triangles $\Omega_{\gamma}=\{(z_1, z_2)\in \mathbb{C}^2: |z_1|^{\gamma}<|z_2|<1\}~(\gamma\in \mathbb{R}^{+})$ investigated by Chakrabarti-Zeytuncu \cite{Chakrabarti-Zeytuncu} and Edholm-McNeal \cite{Edholm-McNeal, Edholm-McNeal1}, Beberok \cite{Beberok} also considered the $L^p$ boundedness of the Bergman projection on the following generalization of the Hartogs triangle
\begin{equation} \label{Bdomain}
\mathcal{H}^{n+1}_{k}:=\{(z,w) \in \mathbb{C}^{n}\times \mathbb{C}: \Vert z \Vert<|w|^{k}<1\},~~k\in \mathbb{Z}^{+},
\end{equation}
where $\Vert \cdot \Vert$ is the Euclidean norm in $\mathbb{C}^n$.
He proved the following result:
\begin{Theorem} $\mathrm{(see \; Beberok \; \cite[Theorem \; 2.1]{Beberok})}$  \label{Btheorem}
Let $p\in (1, \infty)$. Then the Bergman projection is a bounded operator on $L^p(\mathcal{H}_{k}^{n+1})$ if and only if $p\in \big(\frac{2nk+2}{nk+2}, \frac{2nk+2}{nk}\big)$.
\end{Theorem}

Therefore, the restricted range of $p$ is affected by the dimension and the ``camber" of the domain $\mathcal{H}_{k}^{n+1}$.

In recent paper \cite{Chen}, Chen introduced a family of bounded Hartogs domains generalizing the Hartogs triangle as follows.
For $j=1,\cdots, l$, let $\Omega_j$ be a bounded smooth domain in $\mathbb{C}^{k_j}$, $\phi_j:\Omega_j\rightarrow\mathbb{B}^{k_j}$ be a biholomorphic mapping, $m_j=\sum_{s=1}^jk_s$ with $m_0=0$ and $m_l=k_1+\cdots+k_l:=k<n$, $\tilde{z}_j=(z_{m_{j-1}+1},\ldots,z_{m_j})$, and $z=(\tilde{z}_{1}, \cdots, \tilde{z}_{l}, z_{k+1}, \cdots, z_n)\in \mathbb{C}^n$. For $1\leq k<n$, define
\begin{equation}\label{Cdomain}
\mathcal{H}^n_{\{k_j, \phi_j\}}=\{z\in \mathbb{C}^n: \max_{1\leq j \leq l}\Vert\phi_j(\widetilde{z}_j)\Vert<|z_{k+1}|<\cdots<|z_n|<1\},
\end{equation}
when $l=k=1, n=2$, and $\phi_1$ is the identity map, we obtain the classical Hartogs triangle.
The author proved the following result:
\begin{Theorem}  $\mathrm { (see \;  Chen \; \cite[Theorem \; 1.1]{Chen})} $  \label{Ctheorem}
For $1\leq p< \infty$ and $1\leq k<n$, the Bergman projection $P_{\mathcal{H}^n_{\{k_j, \phi_j\}}}$ for $\mathcal{H}^n_{\{k_j, \phi_j\}}$ is bounded on $L^p(\mathcal{H}^n_{\{k_j, \phi_j\}})$ if and only if $p$ is in the range $(\frac{2n}{n+1}, \frac{2n}{n-1})$.
\end{Theorem}

Thus, the $L^p$ boundedness of the Bergman projection on $\mathcal{H}^n_{\{k_j, \phi_j\}}$ is only dependent on the dimension $n$ but not on $\{k_j, \phi_j\}$.

Now, we consider a slightly wider class of non-smooth bounded pseudoconvex domains
which contains the above two domians (\ref{Bdomain}) and (\ref{Cdomain}).
It is defined by
\begin{equation} \label{Wdomain}
\mathcal{H}^{n}_{\{k_j, \phi_j, b\}}=\{z\in \mathbb{C}^n: \max_{1\leq j \leq l}\Vert\phi_j(\widetilde{z}_j)\Vert<|z_{k+1}|^b<\cdots<|z_n|^b<1\},
\end{equation}
where the notations are same as those in (\ref{Cdomain}), and $b\in \mathbb{Z}^+$.
When $k=n-1, ~l=1$, and $\phi_1$ is the identity map, we obtain $\mathcal{H}^{n}_{b}$. There exists a biholomorphism $\Theta:  \mathcal{H}^{n}_{\{k_j, \phi_j, b\}}\rightarrow \mathcal{H}^{n}_{\{k_j, \phi_j\}}$, which is defined by
$$
\Theta(z)=(\phi_1^{-1}(\phi_1(\tilde{z}_{1})z_{k+1}^{1-b}),\cdots, \phi_l^{-1}(\phi_l(\tilde{z}_{l})z_{k+1}^{1-b}), z_{k+1}, \cdots, z_n).
$$
When $b=1$, $\Theta$ is just the identity map.

It is well known that even though two domains are biholomorphic equivalence, the corresponding $L^p$ behavior of the Bergman projection on these two domains may be totally different from each other, let alone the Toeplitz operator constructed by the Bergman projection. Therefore, it would be of interest to consider what will happen for the boundedness of the Bergman projection in this more general setting.

In this paper, inspired by the idea in Chen \cite{Chen} and Khanh-Liu-Thuc \cite{Khanh-Liu-Thuc2}, we mainly focus on the $L^p-L^q$ boundedness of the Toeplitz operator with symbol $K^{-t}(z,z)~(t\geq 0)$ on $\mathcal{H}^n_{\{k_j, \phi_j, b\}}$, where $K^{-t}(z,z):=(K(z,z))^{-t}$, and $K(z, z)$ is the Bergman kernel on diagonal for $\mathcal{H}^n_{\{k_j, \phi_j, b\}}$. We conclude that the parameter $b$ plays an interesting role in the restricted range of $p$ for the boundedness of the Bergman projection on the domain considered in this paper.

\vspace{3mm}
\noindent 1.3. Main results
\begin{Theorem} \label{Wtheorem}
Let $1< p\leq q< \infty$ and $T_{K^{-t}}: L^p(\mathcal{H}_{\{k_j,\phi_j, b\}}^n)\rightarrow L^q(\mathcal{H}_{\{k_j,\phi_j, b\}}^n)$ be the Toeplitz operator with symbol $K^{-t}(z,z)$.
\begin{itemize}
  \vspace{1mm}
  \item[$(1)$] If $q\in \big[\frac{2n+2k(b-1)}{n-1+k(b-1)}, \infty\big)$, then $T_{K^{-t}}$ is unbounded for any $t\in[0, \infty)$.
  \vspace{1mm}
  \item[$(2)$] If $q\in \big(\frac{2(n-1)+2k(b-1)}{n+1+k(b-1)-2/p}, \frac{2n+2k(b-1)}{n-1+k(b-1)}\big)$, then $T_{K^{-t}}$ is bounded if and only if $t\geq\frac{1}{p}-\frac{1}{q}$.
  \vspace{1mm}
  \item[$(3)$] If $q\in \big[p, \frac{2(n-1)+2k(b-1)}{n+1+k(b-1)-2/p}\big]$, then $T_{K^{-t}}$ is bounded if and only if $t>\frac{1}{2p}+\frac{(1-p)}{2p}\frac{n+1+k(b-1)}{n-1+k(b-1)}$.
\end{itemize}
\end{Theorem}

Setting $t=0$ and $p=q$ in Theorem \ref{Wtheorem}, we obtain a sharp range of $p$ for the boundedness of the Bergman projection on $\mathcal{H}_{\{k_j, \phi_j, b\}}^n$ as follows.

\begin{Corollary}\label{Cor1}
Let $1< p<\infty$. Then the Bergman projection from $L^p(\mathcal{H}_{\{k_j,\phi_j, b\}}^n)$ to itself is bounded if and only if $p\in \big(\frac{2n+2k(b-1)}{n+1+k(b-1)}, \frac{2n+2k(b-1)}{n-1+k(b-1)}\big)$.
\end{Corollary}

\begin{Remark}
$(1)$ When $b=1$, Corollary \ref{Cor1} is the known result in Chen \cite[Theorem 1.1]{Chen}
in the case $1<p<\infty$.

$(2)$ When $l=1$, $k=n-1$, and $\phi_1$ is the identity map, Corollary \ref{Cor1} is just the result obtained by Beberok~\cite [Theorem 2.1] {Beberok} for $\mathcal{H}_{b}^n$.
\end{Remark}

It is shown that, influenced by the parameter $b$, the $L^p$-boundedness of the Bergman projections on these domains present an interesting phenomenon. More precisely, comparing with the result obtained by Chen on $\mathcal{H}^n_{\{k_j, \phi_j\}}$ (see Theorem \ref{Ctheorem} mentioned above), the boundedness range of $p$ for the Bergman projection on $\mathcal{H}^n_{\{k_j, \phi_j, b\}}$ is not only dependent on the dimension $n$ but also $k$ and $b$ unless $b=1$, where $k=k_1+\cdots+k_l$ (see Corollary \ref{Cor1}). Moreover, the $L^p$-boundedness range $\big(\frac{2n+2k(b-1)}{n+1+k(b-1)}, \frac{2n+2k(b-1)}{n-1+k(b-1)}\big)$ becomes smaller with the increase of $b$
and approaches $\{2\}$ as $b\rightarrow \infty$.

\vspace{2mm}
Through out this paper, we write $A\lesssim B$ to mean that there exists a constant $C_0>0$ such that $A\leq C_0B$, and use $A\approx B$ for $B\lesssim A\lesssim B$. For the multi-index $\alpha \in \mathbb{Z}^n$, we write $\alpha=(\tilde{\alpha}_1, \cdots, \tilde{\alpha}_l,\alpha_{k+1},\cdots, \alpha_n)\in \mathbb{Z}^{k_1}\times\cdots\times\mathbb{Z}^{k_l}\times{\mathbb{Z}\times\cdots \times\mathbb{Z}}:=\mathbb{Z}^k\times \mathbb{Z}^{n-k}$, where $\tilde{\alpha}_j=(\alpha_{m_{j-1}+1},\ldots,\alpha_{m_j})\in \mathbb{Z}^{k_j}$ with  $m_j=\sum_{s=1}^jk_s$, $m_0=0$ and $m_l=k_1+\cdots+k_l:=k<n$, $j=1, \cdots, l$. Let $|\alpha|=\alpha_1+\cdots+\alpha_n$. In addition, the volume measures mentioned below are all normalized.

\section{Preliminaries}

\setcounter{equation}{0}
\noindent 2.1. The relation between $\mathcal{H}^{n}_{\{k_j, \phi_j, b\}}$ and $\mathcal{H}^{n}_{\{k_j, b\}}$
\vspace{3mm}

When $\phi_j$'s are all identity maps in (\ref{Wdomain}), we denote this special domain by
\begin{equation*}
\mathcal{H}^n_{\{k_j, b\}}=\{z\in \mathbb{C}^n: \max_{1\leq j \leq l}\Vert\widetilde{z}_j\Vert<|z_{k+1}|^b<\cdots<|z_n|^b<1\}.
\end{equation*}
It is easy to verify that $\mathcal{H}^n_{\{k_j, b\}}$ is biholomorphic to $\mathcal{H}^{n}_{\{k_j, \phi_j, b\}}$ by the biholomorphism
\begin{equation}\label{Biholo1}
\Phi(z):=(\phi^{-1}_1(\widetilde{z}_1), \cdots, \phi^{-1}_l(\widetilde{z}_l), z_{k+1}, \cdots, z_n), ~~z \in \mathcal{H}^n_{\{k_j, b\}}.
\end{equation}
The careful reader will note that the biholomorphism (\ref{Biholo1}) is same as the biholomorphism between $\mathcal{H}^n_{\{k_j\}}$ and $\mathcal{H}^n_{\{k_j, \phi_j\}}$, two domains considered by Chen \cite{Chen}.
 Applying Bell's extension Theorem, Chen obtain the equivalence of the $L^p$ boundedness of the Bergman projections on these two domains. Here, this method is also adapted to the special Toeplitz operators on the two domains considered by us.

\begin{Lemma}\label{Lequi}
Let $1<p, q<\infty$, $t\geq 0$, $K_{1}$ and $K$ be the Bergman kernels on diagonal for $\mathcal{H}^n_{\{k_j, b\}}$ and $\mathcal{H}^n_{\{k_j, \phi_j, b\}}$, respectively. Then the following statements are equivalent.
\vspace{1mm}
\begin{itemize}
  \item [(1)] $T_{K_{1}^{-t}}$ is bounded from $L^p\big({\mathcal{H}^n_{\{k_j, b\}}}\big)$ to $L^q\big({\mathcal{H}^n_{\{k_j, b\}}}\big)$.
  \vspace{1mm}
  \item [(2)] $T_{K^{-t}}$ is bounded from $L^p\big({\mathcal{H}^n_{\{k_j, \phi_j, b\}}}\big)$ to $L^q\big({\mathcal{H}^n_{\{k_j, \phi_j, b\}}}\big)$.
\end{itemize}
\end{Lemma}

\noindent Proof.
(1) $\Rightarrow$ (2). We put $\mathcal{H}_1:=\mathcal{H}^n_{\{k_j, b\}}$ and $\mathcal{H}:=\mathcal{H}^n_{\{k_j, \phi_j, b\}}$ for short. From Chen \cite[Section 6]{Chen}, we could find two real numbers $c$ and $d$ such that
\begin{equation} \label{Jacobiestimate}
0<c\leq |\det\Phi'(z)|\leq d, ~~z \in \mathcal{H}_{1}.
\end{equation}
Assume that $T_{K_{1}^{-t}}$ is bounded from $L^p(\mathcal{H}_1)$ to $L^q(\mathcal{H}_1)$. Let $f\in L^p(\mathcal{H})$, $\zeta:=\Phi(z)$, and $\eta:=\Phi(w)$ in the following integral. By the transformation rule for the Bergman kernel under biholomorphism, we have
\begin{eqnarray*}
\Vert T_{K^{-t}}f\Vert^q_{L^q(\mathcal{H})}&=&\int_{{\mathcal{H}}}\Big|\int_{\mathcal{H}}K(\zeta, \eta)f(\eta)K^{-t}(\eta, \eta)dv(\eta)\Big|^qdv(\zeta)\\
&=&\int_{{\mathcal{H}_1}}\Big|T_{K_1^{-t}}\big((f\circ \Phi) \cdot \det \Phi'\cdot|\det \Phi'|^{2t}\big)(z)\Big|^q |\det \Phi'(z)|^{2-q}dv(z)\\
&\lesssim & \Vert T_{K_1^{-t}}\big((f\circ \Phi) \cdot \det \Phi'\cdot|\det \Phi'|^{2t}\big) \Vert^q_{L^q(\mathcal{H}_1)}\\
&\lesssim & \Vert (f\circ \Phi) \cdot \det \Phi'\cdot|\det \Phi'|^{2t} \Vert^p_{L^p(\mathcal{H}_1)}\\
&\lesssim & \Vert f \Vert^p_{L^p(\mathcal{H})}.
\end{eqnarray*}
For the last three lines, we apply the boundedness of $T_{K_1^{-t}}$ and estimate (\ref{Jacobiestimate}).

(2) $\Rightarrow$ (1). Same argument for $\Phi^{-1}$ will show the desired result. This finishes the proof.

\vspace{2mm}
Therefore, it is sufficient to investigate the $L^p-L^q$ boundedness of the Toeplitz operator with symbol
$K_1^{-t} ~(t\geq 0)$ on $\mathcal{H}^{n}_{\{k_j, b\}}$. In the rest of the note, we will focus on the domain
$\mathcal{H}^{n}_{\{k_j, b\}}$.

\vspace{3mm}
\noindent 1.2 The orthogonal basis of $A^2({\mathcal{H}_{\{k_j, b\}}^n})$
\vspace{3mm}

Define a biholomorphic map $\Psi:
\mathcal{H}_{\{k_j, b\}}^n \rightarrow
\Pi^n_{\{k_j\}}$ given by
\begin{equation}\label{Formularbih}
\Psi(\tilde{z}_{1}, \cdots, \tilde{z}_{l}, z_{k+1}, z_n)=\Big(\frac{\tilde{z}_1}{z^b_{k+1}}, \cdots, \frac{\tilde{z}_l}{z^b_{k+1}}, \frac{z_{k+1}}{z_{k+2}}, \cdots, \frac{z_{n-1}}{z_n}, z_n\Big),
\end{equation}
where $\Pi^n_{\{k_j\}}:=\mathbb{B}^{k_1}\times\cdots\times\mathbb{B}^{k_l}\times
\underbrace{\mathbb{D}^*\times\cdots\times\mathbb{D}^*}_{n-k}$.
We denote its inverse by $G$. Then
the determinant of the complex Jacobian of $G$ is given by
\begin{equation}\label{Jacobivalue}
\det G'(\eta)=\prod^n_{j=k+1}\eta_j^{j-1+(b-1)k},~~\eta\in \Pi^n_{\{k_j\}},
\end{equation}
and the Bergman kernel on diagonal for $\mathcal{H}_{\{k_j, b\}}^n$ is
\begin{equation}\label{Bergmankerneldiag}
K_{1}(z, z)=K_{1}(G(\eta), G(\eta))=\frac{1}{\det |G'(\eta)|^2\prod_{j=1}^{l}(1-\Vert\tilde{\eta}_j \Vert^2)^{k_j+1}\prod_{j=k+1}^{n}(1-|\eta_j|^2)^2}.
\end{equation}

\begin{Lemma}\label{Lemmabasis}
Let
$\mathcal{A}:=\{\alpha\in\mathbb{Z}^n: \alpha\in \mathbb{N}^{k}\times\mathbb{Z}^{n-k}\}$.
Then
$$
\mathcal{B}(\mathcal{H}_{\{k_j, b\}}^n):=\big\{z^\alpha: \alpha\in\mathcal{A},~ \sum\nolimits_{j=1}^m\alpha_j+(b-1)\sum\nolimits_{j=1}^k\alpha_j>(1-b)k-m,~ m=k+1, \cdots, n\big\}
$$
is a complete orthogonal basis for $A^2(\mathcal{H}_{\{k_j, b\}}^n)$.
\end{Lemma}

\noindent Proof. We first consider the following operator
$\Gamma_{\Psi}: A^2(\Pi_{\{k_j\}}^n)\rightarrow A^2(\mathcal{H}_{\{k_j, b\}}^n)$
defined by
$$
\Gamma_{\Psi}f:=(f\circ \Psi)\cdot \det\Psi'.
$$
Since $\Psi$ is biholomorphic from $\mathcal{H}_{\{k_j, b\}}^n$ to $\Pi_{\{k_j\}}^n$, it is easy to verify that
$\Gamma_{\Psi}: A^2(\Pi_{\{k_j\}}^n)\rightarrow A^2(\mathcal{H}_{\{k_j, b\}}^n)$ is unitary. Together with the fact that
$
\{e_{\beta}(w):=w^{\beta}=\tilde{w}_1^{\tilde{\beta}_1}\cdots\tilde{w}_l^{\tilde{\beta}_l} w_{k+1}^{\beta_{k+1}}\cdots w_n^{\beta_n},~\beta\in \mathbb{N}^n\}
$
is a complete orthonormal basis in $A^2(\Pi_{\{k_j\}}^n)$. Thus
\begin{equation} \label{Ortho}
\{(e_{\beta}\circ \Psi)(z)\cdot [\det G'(\Psi(z))]^{-1},~\beta\in \mathbb{N}^n\}
\end{equation}
forms a complete orthonormal basis in $A^2(\mathcal{H}_{\{k_j, b\}}^n)$. Substituting (\ref{Formularbih})
and (\ref{Jacobivalue}) into (\ref{Ortho}), we could obtain Lemma \ref{Lemmabasis} after uniting similar terms.

\vspace{3mm}
\noindent 1.3. Other key Lemmas
\vspace{3mm}

We first give the generalised version of Schur's test introduced by Khanh-Liu-Thuc \cite{Khanh-Liu-Thuc1}, which is an important tool of studying the $L^p-L^q$ boundedness for Toeplitz operator with symbol $\psi$.

\begin{Lemma}  $\mathrm { (see \; \cite[Theorem \; 5.1]{Khanh-Liu-Thuc1}) }$    \label{Lemmaschur}
Let $(X, \mu)$, $(Y, \upsilon)$ be measure spaces with $\sigma$-finite, positive measures; let $1<p\leq q<\infty$
and $r\in \mathbb{R}$. Let $K: X\times Y\rightarrow \mathbb{C}$ and $\psi: Y\rightarrow \mathbb{C}$ be measurable functions. Assume that there exist positive measurable functions $h_1, h_2$ on $Y$ and $f$ on $X$
such that
$$
h_1^{-1}h_2\psi\in L^{\infty}(Y, d\upsilon)
$$
and the inequalities
\begin{eqnarray}
&&\int_Y|K(x,y)|^{rp^*}h_1^{p^*}(y)d\upsilon(y)\leq C_1f^{p^*}(x), \label{Schur1}\\
& &\int_X|K(x,y)|^{(1-r)q}f^{q}(x)d\mu(x)\leq C_2h_2^q(y), \label{Schur2}
\end{eqnarray}
hold for almost every $x\in (X, \mu)$ and $y\in (Y, \upsilon)$, where $\frac{1}{p}+\frac{1}{p^*}=1$ and $C_1, C_2$ are positive constants. Then the Toeplitz operator $T_\psi$ associated to the kernel $K$ and the symbol $\psi$ defined by
$$
T_\psi(f)(x):=\int_YK(x,y)f(y)\psi(y)d\upsilon(y)
$$
is bounded from $L^p(Y, \upsilon)$ into $L^q(X, \mu)$. Furthermore,
$$
\Vert T_\psi\Vert_{L^p(Y, \upsilon)\rightarrow L^q(X, \mu)}\leq C_1^{\frac{p-1}{p}}C_2^{\frac{1}{q}}
\Vert h_1^{-1}h_2\psi\Vert_{L^{\infty}(Y, \upsilon)}.
$$
\end{Lemma}

\begin{Lemma}$ \mathrm{ (see \; Herbort \; \cite{Herbort},\;  B{\l}ocki \; \cite{Blocki}) }$ \label{Lgreen}
Let $\Omega$ be a bounded pseudoconvex domain in $\mathbb{C}^n$ and let $s$ be any positive number. Then for any holomorphic function $f$ on $\Omega$ and any $w\in \Omega$,
\begin{equation*}
\int_{\{G(\cdot, w)<-s\}}|f(z)|^2dz\geq e^{-2ns}\frac{|f(w)|^2}{K(w,w)},
\end{equation*}
where $G(\cdot, w)$ is the pluricomplex Green function on $\Omega$ with a pole $w$.
\end{Lemma}

For the definition and properties of the pluricomplex Green function, see Klimek \cite[Chapter 6]{Klimek91}.

Next, using approach as in Khanh-Liu-Thuc \cite[Lemma 3.2]{Khanh-Liu-Thuc2}, we obtain a similar estimate as follows. The difference is that we need to extend their estimates on the unit disc $\mathbb{D}$ to the unit ball during the process of the proof.

\begin{Lemma}\label{Ljacobiestimate}
Fix $w\in \mathcal{H}_{\{k_j, b\}}^n$ and $s\in \mathbb{R}^+$.  Then for any $z \in \{z\in \mathcal{H}_{\{k_j, b\}}^n: G_{\mathcal{H}^n_{\{k_j, b\}}}(z, w)<-s\}$, we have
\begin{equation*}
\frac{K_{1}(z,z)}{K_{1}(w,w)}\approx
\Big|\frac{\det \Psi'(z)}{\det \Psi'(w)}\Big|^2.
\end{equation*}
\end{Lemma}

\noindent Proof. We divide into two steps to prove Lemma \ref{Ljacobiestimate}.

\textbf{Step 1.} Assume that $w \in \mathbb{B}^k$, $\varphi_{w}$ is the automorphism of $\mathbb{B}^k$ taking $0$ to $w$, and $\Vert \varphi_{w}(z)\Vert< e^{-s}$
for all $z \in \mathbb{B}^k$. Then $1-\Vert z \Vert^2\approx 1-\Vert w \Vert^2$. When $k=1$, it will degenerate to the result in \cite[Lemma 3.2]{Khanh-Liu-Thuc2}.

Indeed, employing the properties of the automorphism of the unit ball (see Rudin \cite[Theorem 2.2.2]{Rudin}), we have
$$
\frac{1-\Vert z \Vert^2}{1-\Vert w \Vert^2}=\frac{1-\Vert \varphi_{w}\circ \varphi_{w}(z) \Vert^2}{1-\Vert w \Vert^2}=\frac{1-\Vert \varphi_{w}(z) \Vert^2}{|1-\langle \varphi_{w}(z), w\rangle|^2}.
$$
Hence
$$
\frac{1-\Vert \varphi_{w}(z) \Vert}{1+\Vert \varphi_{w}(z) \Vert}\leq\frac{1-\Vert z \Vert^2}{1-\Vert w \Vert^2}\leq \frac{1+\Vert \varphi_{w}(z) \Vert}{1-\Vert \varphi_{w}(z) \Vert}.
$$
Since $\Vert \varphi_{w}(z)\Vert< e^{-s}$, then there exist positive constants $C_1(s)=\frac{1+e^{-s}}{1-e^{-s}}$, depending only on $s$, such that $C_1^{-1}(s)(1-\Vert w \Vert^2)\lesssim 1-\Vert z \Vert^2 \lesssim C_1(s)(1-\Vert w \Vert^2)$.

\textbf{Step 2.} Since $\mathcal{H}^n_{\{k_j, b\}}$ is holomorphic equivalent to the product domain $\Pi^n_{\{k_j\}}$
via map $\Psi$ (see formula (\ref{Formularbih})), then associating with the biholomorphic invariant and the product-property of the pluricomplex Green function for pseudoconvex domains (see Klimek \cite[Theorem 2.3, Theorem 4.2]{Klimek95}), we obtain that
\begin{eqnarray*}
&&G_{\mathcal{H}^n_{\{k_j, b\}}}(z, w)=G_{\Pi^n_{\{k_j\}}}(\Psi(z), \Psi(w))\\
&=& \max\Big\{G_{\mathbb{B}^{k_1}}\Big(\frac{\tilde{z}_1}{z^b_{k+1}},\frac{\tilde{w}_1}{w^b_{k+1}}\Big), \cdots, G_{\mathbb{D}^{*}}(z_n, w_n)\Big\}
\end{eqnarray*}
for all $z, w\in \mathcal{H}^n_{\{k_j, b\}}$. Then for any  $z \in \{z\in \mathcal{H}^n_{\{k_j, b\}}: G_{\mathcal{H}^n_{\{k_j, b\}}}(z, w)<-s\}$, we have
\begin{eqnarray*}
&&\log\Big\Vert\varphi_{\frac{\tilde{w}_j}{w^b_{k+1}}}\Big(\frac{\tilde{z}_j}{z^b_{k+1}}\Big)\Big\Vert<{-s}, ~j=1, \cdots, l;\\
& &\log\Big\vert\varphi_{\frac{{w}_j}{w_{j+1}}}\Big(\frac{{z}_j}{z_{j+1}}\Big)\Big\vert<{-s},~j=k+1, \cdots, n-1;\\
& &\log\vert\varphi_{w_n}(z_n)\vert<{-s}.
\end{eqnarray*}
Here we use the fact that
$G_{\mathbb{B}^k}(z,w)=\log\Vert \varphi_{w}(z)\Vert$ (see Klimek \cite[Example 2.2]{Klimek95}). By the discussion of Step 1, we get
\begin{eqnarray*}
&&1-\frac{\Vert\tilde{z}_j\Vert^2}{|z_{k+1}|^{2b}}\approx 1-\frac{\Vert\tilde{w}_j\Vert^2}{|w_{k+1}|^{2b}}, ~j=1, \cdots, l;\\
& &1-\Big\vert\frac{{z}_j}{z_{j+1}}\Big\vert^2\approx 1-\Big\vert\frac{{w}_j}{w_{j+1}}\Big\vert^2,~j=k+1, \cdots, n-1;\\
& &1-|z_n|^2\approx 1-|w_n|^2.
\end{eqnarray*}
Thus, applying the explicit formula (\ref{Bergmankerneldiag}) of the Bergman kernel on diagonal for $\mathcal{H}_{\{k_j, b\}}^n$, we could derive the estimate in Lemma \ref{Ljacobiestimate}. The proof is completed.

\begin{Lemma} $\mathrm{ (see\; Khanh-Liu-Thuc \; \cite[Lemma \; 2.5]{Khanh-Liu-Thuc2})\label{Lintestimate1} }$
Let $a\geq 1,~-1<u<0$, and $~c>-2$. Then
$$
\int_{\mathbb{D}}\frac{(1-|w|^2)^u|w|^c}{|1-z\bar{w}|^{2a}}dv(w)\lesssim (1-|z|^2)^{-2a+u+2}, ~~z\in \mathbb{D}.
$$
\end{Lemma}

\begin{Lemma} \label{Lintestimate2}
Let $a\geq 1, ~-1<u<0$. Then
$$
\int_{\mathbb{B}^k}\frac{(1-\Vert w \Vert^2)^u}{|1-\langle z, w \rangle|^{(k+1)a}}dv(w)\lesssim (1-\Vert z\Vert^2)^{u+(k+1)(1-a)}, ~~z\in \mathbb{B}^k.
$$
\end{Lemma}
\noindent Proof. When $a=1$, it is the result in Chen \cite[Lemma 3.2]{Chen}. Since $|1-\langle z, w\rangle|\geq \frac{1}{2}(1-\Vert z \Vert^2)$ for $z, w \in \mathbb{B}^k$, when $a\geq 1$, we have $|1-\langle z, w\rangle|^{(k+1)a} \gtrsim (1-\Vert z \Vert^2)^{(k+1)(a-1)}|1-\langle z, w\rangle|^{k+1}$. Then
\begin{eqnarray*}
&&\int_{\mathbb{B}^k}\frac{(1-\Vert w \Vert^2)^u}{|1-\langle z, w \rangle|^{(k+1)a}}dv(w)\\
&\lesssim& (1-\Vert z\Vert^2)^{(k+1)(1-a)}\int_{\mathbb{B}^k}\frac{(1-\Vert w \Vert^2)^u}{|1-\langle z, w \rangle|^{k+1}}dv(w)\\
&\lesssim& (1-\Vert z\Vert^2)^{(k+1)(1-a)+u},
\end{eqnarray*}
which proves Lemma \ref{Lintestimate2}.

\section{Proof of the main results}

\noindent 3.1. Proof of Theorem \ref{Wtheorem} (1).
\vspace{3mm}

By Lemma \ref{Lequi},  we only need to investigate the $L^p-L^q$ boundedness of the Toeplitz operator with symbol $K_1^{-t}(z, z) ~(t\geq 0)$ on $\mathcal{H}^{n}_{\{k_j, b\}}$.

Here and in the sequel, we set $C_{b,k}:=k(b-1)$. We first prove the following two facts:
\vspace{2mm}
\begin{itemize}
  \item[(a)] If $q\geq\frac{2n+2C_{b,k}}{n-1+C_{b,k}}$, then $z_n^{1-n-C_{b, k}}\notin A^q(\mathcal{H}_{\{k_j, b\}}^n)$.
  \vspace{1mm}
  \item[(b)] If $t\geq 0$ and $z^\beta \in \mathcal{B}(\mathcal{H}_{\{k_j, b\}})$, then $\langle K_1^{-t}(z,z)\bar{z}_n^{n-1+C_{b,k}}, z^{\beta}\rangle_{\mathcal{H}_{\{k_j, b\}}^n}=0$ unless $\beta=(0,\ldots, 0, 1-n-C_{b, k})$.
\end{itemize}
\vspace{1mm}

Indeed, for (a), by Lemma \ref{Lemmabasis}, we learn that $z_n^{1-n-C_{b, k}}\in A^2(\mathcal{H}_{\{k_j, b\}}^n)$, and
\begin{eqnarray*}
&&\int_{\mathcal{H}_{\{k_j, b\}}^n}|z_n^{1-n-C_{b, k}}|^qdv(z)\\
&=&\prod_{j=k+1}^{n-1}\int_{\mathbb{D}^*}|w_{j}|^{2(j-1+C_{b,k})}dv(w_j)\\
& &\cdot\int_{\mathbb{D}^*}|w_n|^{(2-q)(n-1+C_{b,k})}dv(w_n)\\
&<&+\infty
\end{eqnarray*}
if and only if $(2-q)(n-1+C_{b,k})>-2$. This proves (a).

For the second fact (b), let $\beta_{b, k}:=(b-1)(\beta_1+\cdots+\beta_k)$, making the change of variables $z=G(\eta)$, we have
\begin{eqnarray*}
& &\langle K^{-t}(z,z)\bar{z}_n^{n-1+C_{b, k}}, z^{\beta}\rangle_{\mathcal{H}_{\{k_j, b\}}^n}\\
&=&\prod_{j=1}^l\int_{\mathbb{B}^{k_j}}(1-\Vert\tilde{\eta}_j\Vert^2)^{t(k_j+1)} \bar{\tilde{\eta}}_j^{\tilde{\beta}_j}dv({\tilde{\eta}_j})\\
& & \cdot \prod_{j=k+1}^{n-1}\int_{\mathbb{D}^*}\frac{|\eta_j|^
{(2+2t)(j-1+C_{b, k})}}{(1-|\eta_j|^2)^{-2t}}\bar{\eta}_j^{\beta_1+\cdots+\beta_j+\beta_{b, k}}dv(\eta_j)\\
& & \cdot \int_{\mathbb{D}^*}\frac{|\eta_n|^
{(2+2t)(n-1+C_{b, k})}}{(1-|\eta_n|^2)^{-2t}}\bar{\eta}_n^{n-1+C_{b, k}+|\beta|+\beta_{b, k}}dv(\eta_n)\\
&\neq& 0
\end{eqnarray*}
if and only if $\beta=(0, \cdots, 0, 1-n-C_{b, k})$. This proves part (b).

Then by the fact $(b)$, for all $t\geq 0$ and $z^\beta \in \mathcal{B}(\mathcal{H}_{\{k_j, b\}})$, we learn that
\begin{eqnarray*}
&&\langle T_{K_1^{-t}}(\bar{z}_n^{n-1+C_{b, k}}), z^{\beta}\rangle_{\mathcal{H}_{\{k_j, b\}}^n}\\
&=&\langle P_{\mathcal{H}_{\{k_j, b\}}^n}(K_1^{-t}(z,z)\bar{z}_n^{n-1+C_{b, k}}), z^{\beta}
\rangle_{\mathcal{H}_{\{k_j, b\}}^n}\\
&=&\langle K_1^{-t}(z,z)\bar{z}_n^{n-1+C_{b, k}}, z^{\beta}\rangle_{\mathcal{H}_{\{k_j, b\}}^n}\\
&=&0~\text{unless}~\beta=(0,\ldots, 0, 1-n-C_{b, k}).
\end{eqnarray*}
Note that the complete orthogonal basis for $A^2(\mathcal{H}_{\{k_j, b\}}^n)$ in Lemma \ref{Lemmabasis}, then there exists a non-zero constant $C$ such that $T_{K_1^{-t}}(\bar{z}_n^{n-1+C_{b, k}})=Cz_n^{1-n-C_{b, k}}$. We finish the proof of Theorem \ref{Wtheorem} (1) by combining with the above fact (a).

\vspace{3mm}
\noindent 3.2. Proof of Theorem \ref{Wtheorem} (2).
\vspace{3mm}

\textbf{Sufficiency.}  Suppose now that $\frac{2(n-1)+2C_{b,k}}{n+1+C_{b,k}-2/p}<q< \frac{2n+2C_{b,k}}{n-1+C_{b,k}}$ and $t\geq\frac{1}{p}-\frac{1}{q}$.
Denote the Bergman kernel on the the product domain $\Pi_{\{k_j\}}^n$ by $\hat K_1$. Then we consider the following test functions
\begin{eqnarray}\label{Testfunction}
\left \{
\begin{array}{l}
\vspace{1mm}
f(z):=f(G(\eta)):=\rho^{-\lambda}(G(\eta))|\det G'(\eta)|^{-1/p^*};\\
\vspace{1mm}
h_1(w):=h_1({G(\zeta)}):=\rho^{-\lambda}(G(\zeta))\prod_{j=k+1}^n|\zeta_j|^{m_j};\\
\vspace{1mm}
h_2(w):=h_2({G(\zeta)}):=\rho^{-\lambda}(G(\zeta))\hat K_1^{1/p-1/q}(\zeta,\zeta)
          |\det G'(\zeta)|^{-1/p},
\end{array}
\right .
\end{eqnarray}
where $\rho(G(\eta)):=\prod_{j=1}^{l}(1-\Vert\tilde{\eta}_j\Vert^2)\prod_{j=k+1}^{n}(1-|\eta_j|^2)$,
$p^*$ is the conjugate exponent of $p$, and the parameters $\lambda$ and $\{m_j\}_{j=k+1}^n$ satisfy
\begin{eqnarray}\label{Parameter}
\left \{
\begin{array}{l}
\vspace{1mm}
0<\lambda<\min\{\frac{1}{q}, \frac{1}{p^*}\},\\
-\frac{j+1+C_{b, k}}{p^*}<m_j\leq \frac{(j-1+C_{b, k})(q-2p)}{pq}.
\end{array}
\right .
\end{eqnarray}
The existence of $m_j$ follows that
\begin{eqnarray}\label{existence} \nonumber
&&-\frac{j+1+C_{b, k}}{p^*}< \frac{(j-1+C_{b, k})(q-2p)}{pq}\\ \nonumber
& &\Leftrightarrow q>\frac{2(j-1)+2C_{b, k}}{j+1+C_{b, k}-2/p}\\
& &\Leftarrow q>\frac{2(n-1)+2C_{b, k}}{n+1+C_{b, k}-2/p}.
\end{eqnarray}
The last inequality is obvious right according to the sufficient condition of Theorem \ref{Wtheorem} (2).

Let $r=(p^*)^{-1}$ in Lemma \ref{Lemmaschur}. Substituting (\ref{Testfunction}) into (\ref{Schur1}) and (\ref{Schur2}) respectively, and making the change of variables $z=G(\eta)$ and $w=G(\zeta)$, we obtain
\begin{eqnarray}\label{Finterg1} \nonumber
&&\int_{\mathcal{H}_{\{k_j, b\}}^n}|K_1(z, w)|^{rp^*}h_1^{p^*}(w)dv(w)\\ \nonumber
&=&\frac{1}{
|\det G'(\eta)|}\prod_{j=1}^{l}\int_{\mathbb{B}^{k_j}}\frac{(1-\Vert\tilde{\zeta}_j\Vert^2)^{-\lambda p^*}}{|1-\langle \tilde{\eta}_j, \tilde{\zeta}_j\rangle|^{k_j+1}}dv(\tilde{\zeta}_j)\\ \nonumber
& &\cdot\prod_{j=k+1}^{n}\int_{\mathbb{D}^*}\frac{(1-\vert{\zeta}_j\vert^2)^{-\lambda p^*}|\zeta_j|^{m_jp^*+j-1+C_{b, k}}}
{|1-\eta_j\bar{\zeta}_j|^{2}}dv({\zeta}_j)\\
&\lesssim& \rho^{-\lambda p^*}(G(\eta))|\det G'(\eta)|^{-1}=f^{p^*}(z)
\end{eqnarray}
and
\begin{eqnarray} \label{Finterg2} \nonumber
&&\int_{\mathcal{H}_{\{k_j, b\}}^n}|K_1(z, w)|^{(1-r)q}f^{q}(z)dv(z)\\ \nonumber
&=&\frac{1}{
|\det G'(\zeta)|^{q/p}}\prod_{j=1}^{l}\int_{\mathbb{B}^{k_j}}\frac{(1-\Vert\tilde{\eta}_j\Vert^2)^{-\lambda q}}{|1-\langle \tilde{\eta}_j, \tilde{\zeta}_j\rangle|^{(k_j+1)q/p}}dv(\tilde{\eta}_j)\\ \nonumber
& & \cdot \prod_{j=k+1}^{n}\int_{\mathbb{D}^*}\frac{(1-\vert{\eta}_j\vert^2)^{-\lambda q}|\eta_j|^{(2-q)(j-1+C_{b, k})}}
{|1-\eta_j\bar{\zeta}_j|^{2q/p}}dv({\eta}_j)\\
&\lesssim& \hat K_1^{q/p-1}(\zeta,\zeta)\rho^{-\lambda q}(G(\zeta))|\det G'(\zeta)|^{-q/p}=h_2^q(w).
\end{eqnarray}
Here, for the inequalities in formulas (\ref{Finterg1}) and (\ref{Finterg2}), we employ the range of the parameters (\ref{Parameter}) and integral estimates in Lemma \ref{Lintestimate1} and Lemma \ref{Lintestimate2}.

On the other hand, put $w=G(\zeta)$, by the relationship $\hat K_1(\zeta, \zeta)=|\det G'(\zeta)|^2K_1(G(\zeta), G(\zeta))$, we see that
\begin{eqnarray}\label{Foutcome}
&&h_1^{-1}(w)h_2(w)K_1^{-t}(w,w) \nonumber\\
&=&  K_1^{\frac{1}{p}-\frac{1}{q}-t}(G(\zeta),G(\zeta)) \prod\nolimits_{j=k+1}^{n}|\zeta_j|^{(j-1+C_{b, k})(\frac{1}{p}-\frac{2}{q})-m_j}.
\end{eqnarray}
By (\ref{Parameter}), we have $(j-1+C_{b, k})(\frac{1}{p}-\frac{2}{q})-m_j\geq 0$. Again since $t\geq\frac{1}{p}-\frac{1}{q}$,
it follows from formula (\ref{Foutcome}) that $h_1^{-1}h_2K_1^{-t}\in L^{\infty}(\mathcal{H}_{\{k_j, b\}})$. Thus, by Lemma \ref{Lemmaschur}, we obtain that $T_{K_1^{-t}}$ is bounded from $L^p(\mathcal{H}_{\{k_j,b\}}^n)$ to $L^q(\mathcal{H}_{\{k_j,b\}}^n)$. This completes the proof of the sufficiency of Theorem \ref{Wtheorem} (2).

\textbf{Necessity.} Suppose now that $\frac{2(n-1)+2C_{b,k}}{n+1+C_{b,k}-2/p}<q< \frac{2n+2C_{b,k}}{n-1+C_{b,k}}$ and $T_{K_1^{-t}}$ is bounded from $L^p(\mathcal{H}_{\{k_j, b\}}^n)$ to $L^q(\mathcal{H}_{\{k_j, b\}}^n)$.
Let $\mathcal{H}_1:=\mathcal{H}^n_{\{k_j, b\}}$.  We set $g_w(z):=\frac{K_1(z, w)}{\det\Psi'(z)}$ for $z, ~w\in \mathcal{H}_{1}$. Then, similar computation as (\ref{Finterg1}), we have
\begin{eqnarray}\label{Wfunction}
\Vert g_w\Vert_{L^p(\mathcal{H}_1)}\lesssim (K_1(w,w))^{1-\frac{1}{p}}|\det\Psi'(w)|^{\frac{2}{p}-1}.
\end{eqnarray}
Assume that $s\in \mathbb{R}^{+}$. Then, by Lemma \ref{Ljacobiestimate} and Lemma \ref{Lgreen}, we have
\begin{eqnarray*}
&&\int_{\mathcal{H}_1} K_1^{-t}(z,z)|g_w(z)|^2dv(z)\\
&\gtrsim&\int_{\{z\in \mathcal{H}_1, G_{\mathcal{H}_1}(\cdot, w)<-s\}} \left|\frac{K_1(z,w)}{(\det\Psi'(z))^{1+t}}\right|^2
\left|\frac{K_1(z, z)}{(\det\Psi'(z))^2}\right |^{-t}dv(z)\\
&\gtrsim&\left|\frac{K_1(w, w)}{(\det\Psi'(w))^2}\right|^{-t}\int_{\{z\in \mathcal{H}_1, G_{\mathcal{H}_1}(\cdot, w)<-s\}}  \left|\frac{K_1(z,w)}{(\det\Psi'(z))^{1+t}}\right|^2dv(z)\\
&\gtrsim& (K_1(w,w))^{1-t}|\det\Psi'(w)|^{-2}.
\end{eqnarray*}
On the other hand, we also have
\begin{eqnarray*}
&&\int_{\mathcal{H}_1} K_1^{-t}(z,z)|g_w(z)|^2dv(z)\\
&=&\int_{\mathcal{H}_1} K_1^{-t}(z,z)\frac{K_1(w,z)}{\overline{\det\Psi'(z)}}
\frac{K_1(z, w)}{\det\Psi'(z)}dv(z)\\
&=&\int_{\mathcal{H}_1} K_1^{-t}(z,z)\frac{K_1(w,z)}{\overline{\det\Psi'(z)}}
\Big(\int_{\mathcal{H}_1}K_1(z, \eta)\frac{K_1(\eta, w)}{\det\Psi'(\eta)}dv(\eta)\Big)dv(z)\\
&=&\int_{\mathcal{H}_1}\frac{K_1(\eta, w)}{\det\Psi'(\eta)}\Big(\int_{\mathcal{H}_1} K_1(z, \eta)K_1^{-t}(z,z)\frac{K_1(w,z)}{\overline{\det\Psi'(z)}}dv(z)\Big)
dv(\eta)\\
&=&\int_{\mathcal{H}_1}g_w(\eta)\overline{T_{K_1^{-t}}(g_w)(\eta)}dv(\eta)\\
&\lesssim&\Vert g_w\Vert_{L^{q^*}(\mathcal{H}_{1})}
\Vert T_{K_1^{-t}}(g_w)\Vert_{L^q(\mathcal{H}_{1})}\lesssim \Vert g_w\Vert_{L^{q^*}(\mathcal{H}_{1})}
\Vert g_w\Vert_{L^p(\mathcal{H}_{1})}\\
&\lesssim&(K_1(w,w))^{1-\frac{1}{p}+\frac{1}{q}}|\det\Psi'(w)|^{\frac{2}{p}-\frac{2}{q}}.
\end{eqnarray*}
For the last two lines, we apply the H\"{o}lder inequality, the $L^p-L^q$ boundedness of the Toeplitz operator $T_{K_1^{-t}}$, and the estimate (\ref{Wfunction}).
Thus, comparing the above two formulas, we obtain
\begin{equation}\label{Formubijiao}
(K_1(w,w))^{-t+\frac{1}{p}-\frac{1}{q}}|\det\Psi'(w)|^{-2-\frac{2}{p}+\frac{2}{q}}\lesssim \text{constant}.
\end{equation}
Since $|\det\Psi'(w)|>1$ for $w\in\mathcal{H}_{1}$ and $q\geq p > 1$, the second term is bounded. In addition, $K_1(w,w)\rightarrow \infty$ as $w\rightarrow \partial \mathcal{H}_{1}$. Then it follows from (\ref{Formubijiao}) that $\frac{1}{p}-\frac{1}{q}\leq t$. This proves the necessity of Theorem \ref{Wtheorem} (2).

\vspace{3mm}
\noindent 3.3. Proof of Theorem \ref{Wtheorem} (3)
\vspace{3mm}

\textbf{Sufficiency.}  Suppose that $1<p\leq q \leq\frac{2(n-1)+2C_{b,k}}{n+1+C_{b,k}-2/p}$ and $t>\frac{1}{2p}+\frac{(1-p)}{2p}\frac{n+1+C_{b,k}}{n-1+C_{b,k}}$.
In order to prove the sufficiency of Theorem \ref{Wtheorem} (3), we only need to proceed as the proof of the sufficiency of Theorem \ref{Wtheorem} (2). However, we should reset the value of parameters $\{m_j\}_{j=k+1}^n$ by $(j+1+C_{b, k})(1/p-1)<m_j\leq (j-1+C_{b, k})(-1/p+2t)$. Similar as (\ref{existence}), it is easy to derive the existences of $\{m_j\}_{j=k+1}^n$.

On the other hand, we make a deformation of (\ref{Foutcome}) as follows
\begin{eqnarray}\label{Foutcome1}
&&h_1^{-1}(w)h_2(w)K_1^{-t}(w,w) \nonumber\\
&=&\hat{K}_1^{1/p-1/q-t}(\zeta,\zeta)\prod_{j=k+1}^{n}|\zeta_j|^{-m_j+(j-1+C_{b,k})(-1/p+2t)}.
\end{eqnarray}
Since $t>\frac{1}{2p}+\frac{(1-p)}{2p}\frac{n+1+C_{b,k}}{n-1+C_{b,k}}$, it is easy to obtain that $t\geq \frac{1}{p}-\frac{1}{q}$. Combining with the range of $m_j$,
we derive that $h_1^{-1}h_2K_1^{-t}\in L^{\infty}(\mathcal{H}_{\{k_j, b\}})$ from (\ref{Foutcome1}). Thus, by Lemma \ref{Lemmaschur}, we complete the proof of the sufficiency of Theorem \ref{Wtheorem} (3).

\textbf{Necessity.} Suppose that $1<p\leq q \leq\frac{2(n-1)+2C_{b,k}}{n+1+C_{b,k}-2/p}$ and $T_{K_1^{-t}}$ is bounded from $L^p(\mathcal{H}_{\{k_j, b\}}^n)$ to $L^q(\mathcal{H}_{\{k_j, b\}}^n)$. We argue by contradiction. Namely, we assume that $t\leq\frac{1}{2p}+\frac{(1-p)}{2p}\frac{n+1+C_{b,k}}{n-1+C_{b,k}}$. We adopt a family of functions used in Chen \cite {Chen}, which was also applied by Khanh-Liu-Thuc \cite{Khanh-Liu-Thuc2}. Here, in order to hold in our case, we make a little modification of the power of the functions.  Next, for the completeness, we give the details. Define a sequence $\{f_j\}_{j=1}^{\infty}$ by
\begin{eqnarray*}
  f_j(z):=
       \left \{
\begin{array}{l}
h(|z_n|)\bar{z}_n^{n-1+C_{b, k}}; ~~|z_n|\in(a_{j+1},1),\\
\vspace{1mm}
~~~~~~~~~~~~~~~0;~~|z_n|\in(0, a_{j+1}],
\end{array}
      \right .
\end{eqnarray*}
where $a_j:=j^{-j}$ and the function $h:(0,1]\rightarrow (0, \infty)$ is defined by
$$
h(r):=r^{x} ,~~r\in(a_{l+1}, a_l];~~l=1, 2, \cdots,
$$
where $x=\frac{1}{l}-\frac{2}{p}(n+C_{b, k})-(n-1+C_{b, k})$.
A simple calculation shows that $\Vert f_j\Vert_{L^p(\mathcal{H}^n_{\{k_j, b\}})}^p$ is controlled by $\sum_{l=1}^{\infty}l^{-p}$.
Then $f_j \in L^p(\mathcal{H}^n_{\{k_j, b\}})$ for all $p>1$.

On the other hand, we have
\begin{eqnarray} \label{Flast} \nonumber
&&T_{K_{1}^{-t}}(f_j)(G(\eta))=\int_{\mathcal{H}^n_{\{k_j, b\}}}K_{1}(G(\eta), w)K_{1}^{-t}(w, w)f_j(w)dv(w)\\ \nonumber
&=&\int_{\Pi_{\{k_j\}}}K_{1}(G(\eta), G(\zeta))K_{1}^{-t}(G(\zeta), G(\zeta))f_j(G(\zeta))|\det G'(\zeta)|^{2}dv(\zeta)\\
&=&\frac{1}{\det G'(\eta)}\int_{\Pi_{\{k_j\}}}\frac{\hat K_{1}(\eta, \zeta)}{\overline{\det G'(\zeta)}}\hat K_{1}^{-t}(\zeta, \zeta)f_j(G(\zeta))|\det G'(\zeta)|^{2+2t}dv(\zeta).
\end{eqnarray}
Since by the proof of Lemma \ref{Lemmabasis}, $\hat K_{1}(\eta, \zeta)$ could be written as $\sum_{\beta\in \mathbb{N}^n}|c_\beta|^2e_{\beta}(\eta)\overline{e_{\beta}(\zeta)}$. Note that $f_j$ is only dependent on the last variables. Substituting the power series of $\hat K_{1}(\eta, \zeta)$ into (\ref{Flast}). After changing the order of integral and summation, it is easy to obtain that the summation is only work on the index $\beta=(0, \cdots, 0, k+C_{b,k}, \cdots, n-2+C_{b, k}, 0)$. Thus, we learn that
\begin{eqnarray*}\nonumber
&&|T_{K_{1}^{-t}}(f_j)(G(\eta))|\\
&\approx&|\eta_n|^{1-n-C_{b, k}}|\int_{\Pi_{\{k_j\}}}
\overline{\zeta}_{n}^{1-n-C_{b, k}}\hat K_{1}^{-t}(\zeta, \zeta)f_j(G(\zeta))|\det G'(\zeta)|^{2+2t}dv(\zeta)|\\
&\approx&|\eta_n|^{1-n-C_{b, k}}\sum_{l=1}^{j}\int_{a_{l+1}}^{a_l}(1-r)^{2t}r^{x+(2+2t)(n-1+C_{b, k})+1}dr\\
&\gtrsim& \sum_{l=1}^{j}\int_{a_{l+1}}^{a_l}(1-r)^{2t}r^{x+(2+2t)(n-1+C_{b, k})+1}dr.
\end{eqnarray*}
Since $t\leq \frac{1}{2p}+\frac{(1-p)}{2p}\frac{n+1+C_{b,k}}{n-1+C_{b,k}}$, we have $r^{x+(2+2t)(n-1+C_{b, k})+1}\geq r^{1/l-1}$ and $2t<1$. Therefore, we also have $(1-r)^{2t}>1-r$ for any $r\in (0, 1)$.  Then we learn that
$$
\Vert T_{K_1^{-t}}(f_j)\Vert_{L^q({\mathcal{H}^n_{\{k_j, b\}}})}\gtrsim \sum_{l=1}^{j}\int_{a_{l+1}}^{a_l}r^{1/l-1}dr-\sum_{l=1}^{j}\int_{a_{l+1}}^{a_l}r^{1/l}dr.
$$
Since the first term goes to infinity and the second term converges as $j\rightarrow \infty$. So $T_{K_1^{-t}}$ is not bounded, a contradiction. This proves the necessity condition of Theorem
\ref{Wtheorem} (3).
\color{black}

\vskip 5pt

\noindent\textbf{Acknowledgments}\quad  The project is supported
by the National Natural Science Foundation of China (No. 11671306).

\addcontentsline{toc}{section}{References}
\phantomsection
\renewcommand\refname{References}

\clearpage
\end{document}